\documentclass[12pt]{article}
\usepackage{mathrsfs}
\usepackage{amsmath,amssymb}
\newcommand{\be}{\begin{equation}}
\newcommand{\ee}{\end{equation}}

\usepackage{color}

\def\to{\longrightarrow}
\def\mapsto{\longmapsto}

\def\tri{\ \triangle\ }
\def\cwedge{\bigcirc\kern-1.07em\wedge\ }

\newtheorem{thm}{Theorem}[section]

\newtheorem{cor}[thm]{Corollary}

\newtheorem{ex}{Example}
\numberwithin{equation}{section}

\begin{document}
\begin{center}
{\large \bf Four-dimensional almost Hermitian manifolds \\
    with vanishing Tricerri-Vanhecke Bochner curvature tensor}
\end{center}
\footnotetext{\small{\it E-mail addresses}: {\bf prettyfish@skku.edu
}(Y. Euh), {\bf hanada14@skku.edu }(J. Lee), {\bf parkj@skku.edu
}(J. H. Park), {\bf sekigawa@math.sc.niigata-u.ac.jp} (K. Sekigawa)
{\bf ayamada@nagaoka-ct.ac.jp} (A. Yamada).}

\begin{center}
Y. Euh${}^1$, J. Lee$^1$, J. H. Park ${}^{1}$\footnote[2]{Supported
by the Korea Research Foundation Grant funded by the Korean
Government (MOEHRD) KRF-2007-531-C00008.}
 K. Sekigawa$^2$
and A. Yamada$^3$
\end{center}
{\small
\begin{center}
$~~~{}^1$ Department of Mathematics,
    Sungkyunkwan University,
    Suwon, 440-746, KOREA\\
$~~~{}^2$Department of Mathematics,
    Niigata University,
    Niigata 950-2181, JAPAN
\\
$~~~{}^3$Division of General Eduction,
    Nagaoka National College of Technology,
   940-8532, JAPAN
\end{center}

\begin{abstract}
 We study curvature properties
    of four-dimensional almost Hermitian manifolds
    with vanishing Bochner curvature tensor
     as defined by Tricerri and Vanhecke. We give local structure
     theorems for such K\"ahler manifolds, and find out several
     examples related to the theorems.
\end{abstract}


}
\section{Introduction}
\quad\;The Bochner curvature tensor $B$ was defined by Bochner as a
formal analogy of the Weyl conformal curvature tensor \cite{Bo}. The
Bochner K\"ahler manifold which is a K\"ahler manifold with
vanishing Bochner curvature tensor has been studied by Kamishima
\cite{Ka} and Bryant \cite{Br}. Tricerri and Vanhecke \cite{TrVa81}
studied the decomposition of the space of all curvature tensors on a
Hermitian vector space from the view-point of unitary representation
theory and defined a Bochner type conformal curvature tensor $B(R)$
for any almost Hermitian manifold $M=(M,J,g)$. Then tensor field
$B(R)$ is invariant under conformal change of the Riemannian metric
$g$. On one hand, Matsuo \cite{Ma} introduced a generalization of
the Bochner curvature tensor which is called the pseudo-Bochner
curvature tensor on a Hermitian manifold $M=(M,J,g)$ and denoted
with $B_H$, and discussed several curvature properties.

In the present paper, we shall study the curvature properties of
four-dimensional almost Hermitian manifolds with vanishing
Tricerri-Vanhecke Bochner curvature tensor. In the sequel, we shall
call an almost Hermitian manifold with vanishing Tricerri-Vanhecke
Bochner curvature tensor a Tricerri-Vanhecke Bochner flat one, and
also call a four-dimensional almost Hermitian manifold an almost
Hermitian surface.

{The authors are grateful for anonymous referee's helpful comments
concerning this paper.}

\section{Preliminaries}
\quad\;Let $M=(M,J,g)$ be a $2n$-dimensional almost Hermitian
manifold and $\Omega$ the K\"ahler form of $M$ defined by
$\Omega(X,Y)=g(JX,Y)$, for $X$, $Y\in\mathfrak{X}(M)$, where
$\mathfrak{X}(M)$ denotes the Lie algebra of all smooth vector
fields $X$, $Y$ on $M$. We denote by $\nabla$ and $R$ the
Levi-Civita connection and the curvature tensor of $(M,J,g)$ defined
by
\begin{equation}
    R(X,Y)Z
    =[\nabla_X,\nabla_Y]Z-\nabla_{[X,Y]}Z,
\end{equation}

\noindent for $X$, $Y$, $Z\in\mathfrak{X}(M)$. Further, we denote by
$\rho$, $\rho^*$, $\tau$ and $\tau^*$ the Ricci tensor, the Ricci
$*$-tensor, the scalar curvature and the $*$-scalar curvature
defined respectively as:

\begin{equation}
\begin{gathered}
    \rho(X,Y)=\text{tr }(Z\mapsto R(Z,X)Y), \\
    \rho^*(X,Y)=\text{tr }(Z\mapsto R(X,JZ)JY), \\
    \tau=\text{tr }Q, \qquad \tau^*=\text{tr }Q^*
\end{gathered}
\end{equation}
where $Q$ and $Q^*$ are the Ricci operator and the Ricci
$*$-operator defined by $g(QX,Y){}=\rho(X,Y)$ and
$g(Q^*X,Y)=\rho^*(X,Y)$, for $X$, $Y\in\mathfrak{X}(M)$,
respectively. We may easily check that $\rho^*(X,Y)=\rho^*(JY,JX)$
holds for all $X$, $Y\in\mathfrak{X}(M)$, and $\rho^*=\rho$ holds if
$M$ is a K\"ahler manifold. An almost Hermitian manifold M is called
a {\it weakly $*$-Einstein} manifold if $\rho^{*}$ =
$\frac{\tau^{*}}{2n}g$ holds on M and also called a {\it
$*$-Einstein} manifold especially if $\tau^{*}$ is constant. We
denote by $\mathcal{R}$ the curvature operator defined by
\begin{equation}
    g(\mathcal{R}(\iota(x)\wedge\iota(y)),\iota(z)\wedge\iota(w))
    =-g(R(x,y)z,w)
    =-R(x,y,z,w),
\end{equation}
for $x$, $y$, $z$, $w\in T_pM$, $p\in M$, where $\iota$ denotes the
duality : $TM\to\wedge^1M=T^*M$ (the cotangent bundle of $M$). Let
$\{e_i\}$ be an orthonormal basis of $T_pM$ at any point $p\in M$.
In this paper, we shall adopt the following notational convention:
\begin{equation}
\begin{gathered}
    R_{ijkl}=g(R(e_i,e_j)e_k,e_l), \\
    R_{\bar{i}\,jkl}=g(R(Je_i,e_j)e_k,e_l), \\
    \cdots\cdots \\
    R_{\bar{i}\,\bar{j}\,\bar{k}\,\bar{l}}=g(R(Je_i,Je_j)Je_k,Je_l), \\
    \rho_{ij}=\rho(e_i,e_j), \quad
    \cdots, \quad
    \rho_{\bar{i}\,\bar{j}}=\rho(Je_i,Je_j), \\
    \rho^*_{ij}=\rho^*(e_i,e_j), \quad
    \cdots, \quad
    \rho^*_{\bar{i}\,\bar{j}}=\rho^*(Je_i,Je_j), \\
    J_{ij}=g(Je_i,e_j), \qquad
    \nabla_iJ_{jk}=g((\nabla_{e_i}J)e_j,e_k),
\end{gathered}
\end{equation}
and so on, where the Latin indices run over the range
$1,2,\cdots,2n$.

The Bochner curvature tensor $B(R)$ defined by
Tricerri and Vanhecke is stated below:
\begin{equation}
\begin{aligned}
    B(R)
    =&R
        -\frac{1}{4(n+2)(n-2)}\,g\tri\rho
        +\frac{2n-3}{4(n-1)(n-2)}\,g\cwedge\rho \\
    &   -\frac{1}{4(n+2)(n-2)}\,g\tri(\rho J)
        +\frac{1}{4(n-1)(n-2)}\,g\cwedge(\rho J) \\
    &   +\frac{2n^2-5}{4(n+1)(n+2)(n-2)}\,g\tri\rho^*
        -\frac{2n-1}{4(n+1)(n-2)}\,g\cwedge\rho^* \\
    &   +\frac{3}{4(n+1)(n+2)(n-2)}\,g\tri(\rho^*J)
        -\frac{3}{4(n+1)(n-2)}\,g\cwedge(\rho^*J) \\
    &   +\frac{3n\tau-(2n^2-3n+4)\tau^*}{16(n+1)(n+2)(n-1)(n-2)}\,g\tri g
        -\frac{\tau-\tau^*}{8(n-1)(n-2)}\,g\cwedge g
\end{aligned}
\end{equation}
for $n\geqq3$, and
\begin{equation}\label{eq:B(R)4}
\begin{aligned}
    B(R)
    =&R
         + \frac12\,g\cwedge\rho
        + \frac{1}{12}\,\big\{g\triangle\rho^*-g\cwedge\rho^*
           - g\triangle(\rho^*J) + g\cwedge(\rho^*J)\big\} \\
    &   + \frac{3\tau^*-\tau}{96}g\triangle g
        - \frac{\tau+\tau^*}{16}\,g\cwedge g
\end{aligned}
\end{equation}
for $n=2$,
 where for any (0,2)-tensors $a$ and $b$, we set
\begin{equation}
\begin{aligned}
    &(a\cwedge b)(x,y,z,w) \\
    =&a(x,z)b(y,w)-a(x,w)b(y,z)
        +b(x,z)a(y,w)-b(x,w)a(y,z),
\end{aligned}
\end{equation}
\begin{equation}
    \bar{a}(x,y)=a(x,Jy),
\end{equation}
\noindent for $x$, $y$, $z$, $w\in T_pM$, $p\in M$, and we set
\begin{equation}
    a\tri b
    =a\cwedge b+\bar{a}\cwedge\bar{b}+2\bar{a}\otimes\bar{b}+2\bar{b}\otimes\bar{a}.
\end{equation}
\noindent Further, the Weyl curvature tensor is given by
\begin{equation}\label{def:W}
    W
    =R
        +\frac{1}{2n-2}\,g\cwedge\rho
        -\frac{\tau}{2(2n-1)(2n-2)}\,g\cwedge g.
\end{equation}
We denote by $\mathcal{W}$ the Weyl curvature operator.

{ {We note that the Tricerri-Vanhecke
   Bochner curvature tensor $B(R)$ coincides
   with the usual Bochner curvature tensor
   B on
   K\"ahler manifold \cite{TrVa81}.}}
\section{Local structures of Tricerri-Vanhecke
Bochner flat \\ almost Hermitian surfaces} \quad\;In this section,
we shall discuss Tricerri-Vanhecke Bochner flat almost Hermitian
surfaces and give some local structure theorems for these surfaces.
Let $M=(M,J,g)$ be a Tricerri-Vanhecke Bochner flat almost Hermitian
surface. Then, by \eqref{eq:B(R)4}, the curvature tensor $R$ of $M$
can be expressed explicitly by

\begin{equation}
\begin{aligned}
\label{eq:R_B(R)=0}
    &R(X,Y,Z,W) \\
    =&\frac12\Big\{g(X,W)\rho(Y,Z)+g(Y,Z)\rho(X,W)  \\
    &       \qquad -g(X,Z)\rho(Y,W)-g(Y,W)\rho(X,Z)\Big\} \\
    &   +\frac{1}{12}\Big\{2g(X,JY)\Big(\rho^*(W,JZ)-\rho^*(JZ,W)\Big)  \\
    &       \qquad\quad +2g(Z,JW)\Big(\rho^*(Y,JX)-\rho^*(JX,Y)\Big)  \\
    &       \qquad\quad +g(X,JZ)\Big(\rho^*(W,JY)-\rho^*(JY,W)\Big) \\
    &       \qquad\quad +g(Y,JW)\Big(\rho^*(Z,JX)-\rho^*(JX,Z)\Big) \\
    &       \qquad\quad +g(X,JW)\Big(\rho^*(Y,JZ)-\rho^*(JZ,Y)\Big)  \\
    &       \qquad\quad +g(Y,JZ)\Big(\rho^*(X,JW)-\rho^*(JW,X)\Big)\Big\}  \\
    &   +\frac{3\tau^*-\tau}{48}\Big\{g(X,W)g(Y,Z)-g(X,Z)g(Y,W)  \\
    &       \qquad\qquad\quad -2g(X,JY)g(Z,JW)-g(X,JZ)g(Y,JW)  \\
    &       \qquad\qquad\quad +g(Y,JZ)g(X,JW)\Big\} \\
    &   -\frac{\tau+\tau^*}{8}\Big\{g(X,W)g(Y,Z)-g(X,Z)g(Y,W)\Big\}
\end{aligned}
\end{equation}
 for $X$, $Y$, $Z$, $W\in\mathfrak{X}(M)$. On one hand,
from \eqref{def:W}, the Weyl curvature tensor $W$ is given by
\begin{equation}
\begin{aligned}\label{def:W4}
    &W(X,Y,Z,W) \\
    =&R(X,Y,Z,W)  \\
    &   -\frac12\Big\{g(X,W)\rho(Y,Z)+g(Y,Z)\rho(X,W)  \\
    &       \qquad -g(X,Z)\rho(Y,W)-g(Y,W)\rho(X,Z)\Big\}  \\
    &   +\frac{\tau}{6}\{g(X,W)g(Y,Z)-g(X,Z)g(Y,W)\}
\end{aligned}
\end{equation}
for $X$, $Y$, $Z$, $W\in\mathfrak{X}(M)$. From \eqref{eq:R_B(R)=0}
and \eqref{def:W4}, the Weyl curvature tensor $W$ is also expressed
by
\begin{equation}
\begin{aligned}\label{eq:W_B(R)=0}
    &W(X,Y,Z,W) \\
    =&\frac{\tau-3\tau^*}{24}\Big\{g(X,W)g(Y,Z)-g(X,Z)g(Y,W)\Big\}  \\
    &   +\frac{1}{12}\Big\{2g(X,JY)\Big(\rho^*(W,JZ)-\rho^*(JZ,W)\Big)  \\
    &       \qquad\quad +2g(Z,JW)\Big(\rho^*(Y,JX)-\rho^*(JX,Y)\Big) \\
    &       \qquad\quad +g(X,JZ)\Big(\rho^*(W,JY)-\rho^*(JY,W)\Big) \\
    &       \qquad\quad +g(Y,JW)\Big(\rho^*(Z,JX)-\rho^*(JX,Z)\Big) \\
    &       \qquad\quad +g(X,JW)\Big(\rho^*(Y,JZ)-\rho^*(JZ,Y)\Big)  \\
    &       \qquad\quad +g(Y,JZ)\Big(\rho^*(X,JW)-\rho^*(JW,X)\Big)\Big\} \\
    &   +\frac{3\tau^*-\tau}{48}\Big\{g(X,W)g(Y,Z)-g(X,Z)g(Y,W) \\
    &       \qquad\qquad\quad -2g(X,JY)g(Z,JW)-g(X,JZ)g(Y,JW) \\
    &       \qquad\qquad\quad +g(Y,JZ)g(X,JW)\Big\}
\end{aligned}
\end{equation}
for $X$, $Y$, $Z$, $W\in\mathfrak{X}(M)$. First, from
\eqref{eq:R_B(R)=0}, by direct calculation, we have the following
{\ theorem}.
\begin{thm}\label{th:curv_id}
    Let $M=(M,J,g)$ be a Tricerri-Vanhecke Bochner flat
    almost Hermitian surface.
    Then, the curvature tensor $R$ satisfies the following curvature identity
\begin{equation}\label{eq:curv_id}
\begin{aligned}
    &R(X,Y,Z,W)-R(JX,JY,Z,W) \\
    &   \quad -R(X,Y,JZ,JW)+R(JX,JY,JZ,JW) \\
    =&R(X,JY,Z,JW)+R(X,JY,JZ,W) \\
    &   \quad +R(JX,Y,JZ,W)+R(JX,Y,Z,JW)
\end{aligned}
\end{equation}
    for $X$, $Y$, $Z$, $W\in\mathfrak{X}(M)$.
\end{thm}

\noindent{\bf Proof }
 From \eqref{eq:R_B(R)=0}, the left-hand side of
\eqref{eq:curv_id} is
\begin{equation}
\begin{aligned}\label{eq:curv_id_l}
    &R(X,Y,Z,W)-R(JX,JY,Z,W) \\
    &   \quad -R(X,Y,JZ,JW)+R(JX,JY,JZ,JW)  \\
    =&\frac12\Big\{g(X,W)\rho(Y,Z)+g(Y,Z)\rho(X,W)\\
    &       \qquad -g(X,Z)\rho(Y,W)-g(Y,W)\rho(X,Z) \\
    &       \qquad -g(X,JW)\rho(Y,JZ)-g(Y,JZ)\rho(X,JW) \\
    &       \qquad +g(X,JZ)\rho(Y,JW)+g(Y,JW)\rho(X,JZ)  \\
    &       \qquad +g(X,W)\rho(JY,JZ)+g(Y,Z)\rho(JX,JW)\\
    &       \qquad -g(X,Z)\rho(JY,JW)-g(Y,W)\rho(JX,JZ)  \\
    &       \qquad +g(X,JW)\rho(JY,Z)+g(Y,JZ)\rho(JX,W)\\
    &       \qquad -g(X,JZ)\rho(JY,W)-g(Y,JW)\rho(JX,Z)\Big\} \\
    &   -\frac{\tau+\tau^*}{4}\Big\{g(X,W)g(Y,Z)-g(X,Z)g(Y,W)  \\
    &       \qquad\quad -g(X,JW)g(Y,JZ)+g(X,JZ)g(Y,JW) \Big\}
\end{aligned}
\end{equation}
and the right-hand side of \eqref{eq:curv_id} is

\begin{equation}
\begin{aligned}\label{eq:curv_id_r}
    &R(X,JY,Z,JW)+R(X,JY,JZ,W) \\
    &   \quad +R(JX,Y,JZ,W)+R(JX,Y,Z,JW)  \\
    =&\frac12\Big\{g(X,W)\rho(Y,Z)+g(Y,Z)\rho(X,W) \\
    &       \qquad -g(X,Z)\rho(Y,W)-g(Y,W)\rho(X,Z)  \\
    &       \qquad -g(X,JW)\rho(Y,JZ)-g(Y,JZ)\rho(X,JW)  \\
    &       \qquad +g(X,JZ)\rho(Y,JW)+g(Y,JW)\rho(X,JZ)  \\
    &       \qquad +g(X,W)\rho(JY,JZ)+g(Y,Z)\rho(JX,JW)  \\
    &       \qquad -g(X,Z)\rho(JY,JW)-g(Y,W)\rho(JX,JZ)  \\
    &       \qquad +g(X,JW)\rho(JY,Z)+g(Y,JZ)\rho(JX,W)  \\
    &       \qquad -g(X,JZ)\rho(JY,W)-g(Y,JW)\rho(JX,Z)\Big\}\\
    &   -\frac{\tau+\tau^*}{4}\Big\{g(X,W)g(Y,Z)-g(X,Z)g(Y,W)  \\
    &       \qquad\quad -g(X,JW)g(Y,JZ)+g(X,JZ)g(Y,JW)  \Big\}.
\end{aligned}
\end{equation}
 From \eqref{eq:curv_id_l} and \eqref{eq:curv_id_r}, we may see
that the curvature identity \eqref{eq:curv_id} holds.
\hfill$\square$\medskip


\noindent{\bf Remark} {\rm    It is known that the curvature tensor
of any Hermitian manifold
    satisfies the curvature identity \eqref{eq:curv_id}
    in the above Theorem \ref{th:curv_id} \cite{Gr}.
    However, the converse is not true in general.
    In fact, Tricerri and Vanhecke { \cite{TrVa77}} gave an example
    of a locally flat almost Hermitian surface
    which is not Hermitian.}

Now, let $\{e_i\}=\{e_1,e_2=Je_1,e_3,e_4=Je_3\}$ be a unitary basis
(resp. any local unitary frame field) of $T_pM$ { for any $p\in M$},
and $\{e^i\}$ be the dual basis (resp. local dual unitary frame
field) of $\{e_i\}$. The space $\wedge^2_pM$ of all 2-forms on $M$
is decomposed by
\begin{equation}
    \wedge^2_pM
    =\wedge^2_+\oplus\wedge^2_-,
\end{equation}
and these subspaces are spanned respectively by
\begin{equation}\label{eq:wedge2+-}
    \wedge^2_+
    =\text{span}\left\{\Omega_0,\Phi,J\Phi\right\}, \qquad
    \wedge^2_-
    =\text{span}\left\{\Psi_1,\Psi_2,\Psi_3\right\},
\end{equation}
where
\begin{equation}\label{def:base_wedge2}
\begin{gathered}
    \Omega_0
    =\frac{1}{\sqrt{2}}\Omega=\frac{1}{\sqrt{2}}(e^1\wedge e^2+e^3\wedge e^4), \\
    \Phi
    =\frac{1}{\sqrt{2}}(e^1\wedge e^3-e^2\wedge e^4), \qquad
    J\Phi
    =\frac{1}{\sqrt{2}}(e^1\wedge e^4+e^2\wedge e^3), \\
    \Psi_1
    =\frac{1}{\sqrt{2}}(e^1\wedge e^2-e^3\wedge e^4), \\
    \Psi_2
    =\frac{1}{\sqrt{2}}(e^1\wedge e^3+e^2\wedge e^4), \qquad
    \Psi_3
    =\frac{1}{\sqrt{2}}(e^1\wedge e^4-e^2\wedge e^3).
\end{gathered}
\end{equation}
Then, from \eqref{eq:W_B(R)=0} and \eqref{def:base_wedge2}, we have
\begin{equation}
\begin{aligned}\label{eq:W_op}
    \mathcal{W}(\Omega_0)
    =&\frac{3\tau^*-\tau}{12}\Omega_0
        -\frac12(\rho^*_{14}-\rho^*_{41})\Phi
        +\frac12(\rho^*_{13}-\rho^*_{31})J\Phi, \\
    \mathcal{W}(\Phi)
    =&-\frac12(\rho^*_{14}-\rho^*_{41})\Omega_0
        -\frac{3\tau^*-\tau}{24}\Phi,  \\
    \mathcal{W}(J\Phi)
    =&\frac12(\rho^*_{13}-\rho^*_{31})\Omega_0
        -\frac{3\tau^*-\tau}{24}J\Phi,  \\
    \mathcal{W}(\Psi_i)
    =&0, \qquad(i=1,2,3),
\end{aligned}
\end{equation}
 where $\mathcal{W}$ is the Weyl curvature operator. Thus, by
\eqref{eq:wedge2+-} and \eqref{eq:W_op}, we have the following
theorems.

\begin{thm}\label{th:self-dual}
    Let $M=(M,J,g)$ be a Tricerri-Vanhecke Bochner flat
    almost Hermitian surface.
    Then, $M$ is self-dual.
\end{thm}

\begin{thm}
    Let $M=(M,J,g)$ be a Tricerri-Vanhecke Bochner flat
    almost Hermitian surface. Then, $M$ is anti-self-dual if and only if
    $\rho^*$ is symmetric and $3\tau^*-\tau=0$ holds on $M$.
\end{thm}

\begin{cor}\label{cor:conformally_flat}
    Let $M=(M,J,g)$ be a Tricerri-Vanhecke Bochner flat
    almost Hermitian surface.
    Then, $M$ is conformally flat if and only if
    $\rho^*$ is symmetric and $3\tau^*-\tau=0$ holds on $M$.
\end{cor}

Below are two examples of  conformally flat, Tricerri-Vanhecke
Bochner flat almost Hermitian surfaces.

\begin{ex}
{\rm Let $M=\mathbb{R}^4_+ =\{\,(x_1,x_2,x_3,x_4)\in\mathbb{R}^4
\,|\, x_4>0,\, x_1, x_2, x_3 \in\mathbb{R}\,\}$ and
$\{e_1,e_2,e_3,e_4\}$ be the global frame field on $M$ defined by
\begin{equation}
    e_1=x_4\frac{\partial}{\partial x_1}, \qquad
    e_2=x_4\frac{\partial}{\partial x_2}, \qquad
    e_3=x_4\frac{\partial}{\partial x_3}, \qquad
    e_4=x_4\frac{\partial}{\partial x_4}.
\end{equation}
Further, we define almost Hermitian structure $(J,g)$ on $M$ as
follows:
\begin{equation}
    J: \quad
    e_1\mapsto e_2, \qquad
    e_2\mapsto -e_1, \qquad
    e_3\mapsto e_4, \qquad
    e_4\mapsto -e_3,
\end{equation}
and
\begin{equation}
    g(e_i,e_j)=\delta_{ij}.
\end{equation}
Then, we may easily check that $(M,J,g)$  is a Hermitian surface of
constant sectional curvature $-1$ (and hence, conformally flat,
Tricerri-Vanhecke Bochner flat Hermitian surface by virtue of
\eqref{eq:R_B(R)=0}).}
\end{ex}

\begin{ex}{\rm {\cite{Ta} }}\label{ex:Tanno}
{\rm Let $M_1=(M_1(K),J_1,g_1)$, $M_2=(M_2(-K),J_2,g_2)$ be oriented
surfaces with constant Gaussian curvatures $K$ and $-K$ ($K>0)$
respectively, and $(M,J,g)=(M_1\times M_2,J_1\times J_2,g_1\times
g_2)$ be the direct product of $M_1$ and $M_2$.}
\end{ex}
\noindent We may immediately observe that a complex space form is a
typical example of Tricerri-Vanhecke Bochner flat K\"ahler manifold
and the above Example \ref{ex:Tanno} is such an example. Now,
concerning the Example \ref{ex:Tanno}, we have the following {
theorem.}

\begin{thm}\label{th: K or -K}
    Let $M=(M,J,g)$ be a Tricerri-Vanhecke Bochner flat
    K\"ahler surface.
    If the scalar curvature $\tau$ of $M$ is constant,
    then $M$ is locally a complex space form of complex dimension
    $2$,
    or locally a product of two oriented surfaces
    of different constant Gaussian curvatures $K$ and $-K$ \rm{(}$K\neq0$\rm{)}.
\end{thm}

{\noindent\bf Proof }
 Let $\lambda$, $\mu$
($\lambda\geqq\mu$) be the eigenvalues of the Ricci transformation
$Q$ at each point of $M$. Then, we may easily observe that
$\lambda+\mu=\dfrac{\tau}{2}$ and the eigenvalues $\lambda$, $\mu$
give rise to continuous functions on $M$. Now, we set $M_0=\{\,p\in
M \,|\, \lambda>\mu \text{ at } p\,\}$. Then, $M_0$ is an open set
(possibly, empty set) of $M$.

First, we assume that $M_0$ is empty. Then, we see that $M$ is
Einstein, and hence, by \eqref{eq:R_B(R)=0}, $M$ is locally a
complex space form of complex dimension 2 of constant holomorphic
sectional curvature $\dfrac{\tau}{6}$.

Next, we assume that $M_0$ is not empty. Then, we may define two
smooth $J$-invariant distributions $D_\lambda$ and $D_\mu$ on $M_0$
corresponding to the eigenvalues $\lambda$ and $\mu$ of the Ricci
transformation $Q$. Let $U$ be an any component of $M_0$ and
$\{e_i\}=\{e_1, e_2=Je_1, e_3, e_4=Je_3\}$ be any local unitary
frame field in $U$ such that $Qe_1=\lambda e_1$ ($Qe_2=\lambda
e_2$), $Qe_3=\mu e_3$ ($Qe_4=\mu e_4$). We set
\begin{equation}
    \nabla_{e_i}e_j=\sum_k \Gamma_{ijk}e_k
    \qquad (i,j=1,2,3,4).
\end{equation}
Then, since $M$ is K\"ahler, we get
\begin{equation}\label{eq:Gamma}
    \Gamma_{ijk}=-\Gamma_{ikj}, \qquad
    \Gamma_{i\,\bar{j}\,\bar{k}}=\Gamma_{ijk} \qquad
    (i,j,k=1,2,3,4).
\end{equation}
 On one hand, from \eqref{eq:R_B(R)=0}, since $\tau$ is constant,
we have
\begin{equation}
    \sum_i(\nabla_{e_i}R)(X,Y,Z,e_i)
    =\frac12\left\{(\nabla_X\rho)(Y,Z)-(\nabla_Y\rho)(X,Z)\right\}
\end{equation}
and hence, taking account of the second Bianchi identity,
\begin{equation}\label{eq:nabla_rho}
    (\nabla_X\rho)(Y,Z)-(\nabla_Y\rho)(X,Z)=0
\end{equation}
for any $X$, $Y$, $Z\in\mathfrak{X}(U)$. Thus, by setting
$(X,Y,Z)=(e_1,e_2,e_2)$, $(e_1,e_3,e_3)$, $(e_1,e_4,e_4)$ in
\eqref{eq:nabla_rho}, from \eqref{eq:Gamma}, we have respectively
\begin{gather}
    e_1\lambda=0, \label{eq:122} \\
    e_1\mu+\Gamma_{342}(\lambda-\mu)=0, \label{eq:133} \\
    e_1\mu-\Gamma_{432}(\lambda-\mu)=0. \label{eq:144}
\end{gather}
From \eqref{eq:122} and the hypothesis ($\tau=2(\lambda+\mu)$ is
constant), we have $e_1\mu=0$. Thus, by \eqref{eq:133} and
\eqref{eq:144}, we have $\Gamma_{342}=\Gamma_{432}=0$. Similarly, we
have $e_a\lambda=e_a\mu=0$ ($a=2,3,4$), and
$\Gamma_{341}=\Gamma_{431}=\Gamma_{142}=\Gamma_{241}=\Gamma_{132}=\Gamma_{231}=0$.
Thus, we see that the distributions $D_\lambda$ and $D_\mu$ are both
parallel ones on each $U$. Therefore, $M_0$ is locally { a} product
of two integral manifold{s} with respect to the distributions
$D_\lambda$ and $D_\mu$. From \eqref{eq:R_B(R)=0},
$0=R_{1313}=-\dfrac{\tau}{24}$ and hence $\tau=0$. Since $\lambda$
and $\mu$ are both constant, by setting $\lambda=K$ and $\mu=-K$
from taking account of $\lambda+\mu=\dfrac{\tau}{2}=0$, we see that
$M_0=M$ and $M$ is  locally a product of two oriented surfaces of
different constant Gaussian curvatures $K$ and $-K$
\rm{(}$K\neq0$\rm{)}. \hfill$\square$\medskip

{{ We note that Bochner flat K\"ahler manifold with constant scalar
curvature is locally symmetric in any dimension \cite{MT}.}}

The following example illustrates the above Theorem \ref{th: K or
-K}. Namely, there exists a Tricerri-Vanhecke Bochner flat almost
K\"ahler surface with constant scalar curvature which is not a
K\"ahler one.
\begin{ex}
{\rm We set $(M, g)=\mathbb{H}^3(-1)\times \mathbb{R}$, where
$\mathbb{H}^3(-1)$ is a 3-dimensional real hyperbolic space of
constant sectional curvature $-1$ and $\mathbb{R}$ is a real line.
Let
\begin{equation*}
    e_1=x_1\frac{\partial}{\partial x_1}, \qquad
    e_2=x_1\frac{\partial}{\partial x_2}, \qquad
    e_3=x_1\frac{\partial}{\partial x_3}, \qquad
    e_4=\frac{\partial}{\partial x_4}.
\end{equation*}
on $M=\mathbb{R}^4_+ =\mathbb{R}^3_+\times \mathbb{R}
=\{\,(x_1,x_2,x_3,x_4)\in\mathbb{R}^4 \,|\, x_1>0 \}$ and define an
almost Hermitian structure $(J, g)$ on $M$ by
$g(e_i,e_j)=\delta_{ij}$ and $Je_i=\sum_{j=1}^{4}J_{ij}e_j$, where

\begin{equation*}
(J_{ij}) =\begin{pmatrix}
  0             & \text{cos} x_4      & \text{sin}x_4    & 0      \\
  -\text{cos} x_4      & 0            & 0         & -\text{sin}x_4\\
  -\text{sin}x_4       & 0            & 0         & \text{cos}x_4 \\
  0             & \text{sin}x_4       & -\text{cos} x_4  & 0
 \end{pmatrix}.
\end{equation*}

\noindent We denote by $\{e^i\}_{i=1, \cdots , 4}$ the dual basis of
$\{e_i\}$. Then the K\"ahler form $\Omega$ is given by

\begin{equation}
\begin{aligned}
\Omega
    =J_{12}e^1 \wedge e^2+J_{13}e^1 \wedge e^3+J_{14}e^1 \wedge
    e^4\\
    + J_{23}e^2 \wedge e^3+J_{24}e^2 \wedge e^4+J_{34}e^3 \wedge e^4\\
    \qquad =\frac{1}{x_1^2}\text{cos}x_4 dx_1 \wedge dx_2 +
    \frac{1}{x_1^2} \text{sin}x_4 dx_1 \wedge dx_3\\
    -\frac{1}{x_1}\text{sin}x_4dx_2 \wedge dx_4 +
    \frac{1}{x_1}\text{cos}x_4dx_3 \wedge dx_4
\end{aligned}
\end{equation}
Thus, we have $d\Omega=0$, and hence $(M, J, g)$ is an almost
K\"ahler manifold.}
\end{ex}
We may easily check that Example 3 is a locally symmetric,
conformally flat, Tricerri-Vanhecke Bochner flat, non-K\"ahler,
almost K\"ahler surface with constant scalar curvature $\tau=-6$ and
constant $*$-scalar curvature $\tau^*=-2$.

\section{Compact Tricerri-Vanhecke
Bochner flat almost Hermitian surfaces}
 \quad\;
Let $M=(M,J,g)$ be a compact Tricerri-Vanhecke Bochner flat almost
Hermitian surface. From \eqref{eq:W_op}, we have
\begin{equation}\label{eq:|W_op|}
\begin{aligned}
    \|\mathcal{W}_+\|^2
    =&\frac{(3\tau^*-\tau)^2}{96}
        +\frac12\left\{(\rho^*_{13}-\rho^*_{31})^2
            +(\rho^*_{14}-\rho^*_{41})^2\right\}, \\
    \|\mathcal{W}_-\|^2
    =&0.
\end{aligned}
\end{equation}
We set
\begin{equation}\label{def:G}
    G
    =\sum_{i,j}\left(\rho^*_{ij}-\rho^*_{ji}\right)^2
    =4\left\{\left(\rho^*_{13}-\rho^*_{31}\right)^2
        +\left(\rho^*_{14}-\rho^*_{41}\right)^2\right\}.
\end{equation}
From \eqref{eq:|W_op|}, taking account of \eqref{def:G}, the first
Pontrjagin number is given by

\begin{equation}\label{eq:p_1(M)}
\begin{aligned}
    p_1(M)
    =&\frac{1}{4\pi^2}\int_M\left\{\|\mathcal{W}_+\|^2-\|\mathcal{W}_-\|^2\right\}dv \\
    =&\frac{1}{4\pi^2}\int_M\left\{\frac{(3\tau^*-\tau)^2}{96}+\frac{G}{8}\right\}dv \\
    =&\frac{1}{32\pi^2}\int_M\left\{\frac{(3\tau^*-\tau)^2}{12}+G\right\}dv.
\end{aligned}
\end{equation}
From \eqref{def:W4} and \eqref{eq:W_B(R)=0}, taking account of
\eqref{def:G}, we have

\begin{equation}\label{eq:|R|^2}
\begin{aligned}
    \|R\|^2
    =&\frac{1}{24}\left\{9(\tau^*)^2-6\tau\tau^*-31\tau^2\right\}
        +\sum_{i<j}(\rho_{ii}+\rho_{jj})^2
        +4\sum_{i<j}\rho_{ij}^2 \\
    &   +2\left\{(\rho^*_{13}-\rho^*_{31})^2
            +(\rho^*_{14}-\rho^*_{41})^2\right\} \\
    =&\frac{1}{24}\left\{9(\tau^*)^2-6\tau\tau^*-31\tau^2\right\}
        +2\|\rho\|^2
        +\tau^2 \\
    &   +2\left\{(\rho^*_{13}-\rho^*_{31})^2
            +(\rho^*_{14}-\rho^*_{41})^2\right\} \\
    =&\frac{1}{24}(3\tau^*-\tau)^2-\frac{4}{3}\tau^2
        +2\left(\Big\|\rho-\frac{\tau}{4}g\Big\|^2+\frac{\tau^2}{4}\right)
        +\tau^2
        +\frac12G \\
    =&\frac{1}{24}(3\tau^*-\tau)^2
        +2\Big\|\rho-\frac{\tau}{4}g\Big\|^2
        +\frac{\tau^2}{6}
        +\frac12G.
\end{aligned}
\end{equation}

\noindent From \eqref{eq:|R|^2}, the Euler number is given by

\begin{equation}\label{eq:chi(M)}
\begin{aligned}
    \chi(M)
    =&\frac{1}{32\pi^2}\int_M\left\{\|R\|^2-4\|\rho\|^2+\tau^2\right\}dv \\
    =&\frac{1}{32\pi^2}\int_M\left\{\frac{(3\tau^*-\tau)^2}{24}
        -2\Big\|\rho-\frac{\tau}{4}g\Big\|^2
        +\frac{\tau^2}{6}
        +\frac12G\right\}dv.
\end{aligned}
\end{equation}

\noindent From \eqref{eq:p_1(M)} and \eqref{eq:chi(M)}, by Wu's
theorem \cite{Wu}, the first Chern number is given by

\begin{equation}\label{eq:c_1(M)^2}
\begin{aligned}
    c_1(M)^2
    =&p_1(M)+2\chi(M) \\
    =&\frac{1}{32\pi^2}\int_M\left\{\frac{(3\tau^*-\tau)^2}{6}
        -4\Big\|\rho-\frac{\tau}{4}g\Big\|^2
        +\frac{\tau^2}{3}
        +2G\right\}dv.
\end{aligned}
\end{equation}

\noindent  From \eqref{eq:p_1(M)} and \eqref{eq:c_1(M)^2}, we have
the following results.

\begin{thm}
    Let $M=(M,J,g)$ be  a compact Tricerri-Vanhecke Bochner flat
    almost Hermitian Einstein surface.
    If the first Pontrjagin { number} $p_1(M)$ of $M$ vanishes,
    then $M$ is a space of constant sectional curvature $\dfrac{\tau}{12}$
    \rm{(}$\tau\leqq0$\rm{)}.
\end{thm}

{\noindent\bf Proof }
 From \eqref{eq:p_1(M)}, we have
$3\tau^*-\tau=0$ and $G=0$ which implies $\rho^*$ is symmetric.
Thus, from Corollary \ref{cor:conformally_flat}, $M$ is conformally
flat and hence, $M$ is a space of constant sectional curvature
$\dfrac{\tau}{12}$ since $M$ is Einstein. It is well-known that a
four-dimensional sphere $S^4$ can not admit an almost complex
structure. Therefore, it follows that $\tau\leqq0$.
\hfill$\square$\medskip


\begin{thm}
    Let $M=(M,J,g)$  be a compact Tricerri-Vanhecke Bochner flat
    almost Hermitian Einstein surface.
    If the first Chern number $c_1(M)^2$ of $M$ vanishes,
    then $M$ is locally flat.
\end{thm}

{\noindent\bf Proof }
 Since $M$ is Einstein, from \eqref{eq:c_1(M)^2}, we
have
\begin{equation}
    c_1(M)^2
    =\frac{1}{32\pi^2}\int_M\left\{\frac{(3\tau^*-\tau)^2}{6}
        +\frac{\tau^2}{3}
        +2G\right\}dv.
\end{equation}
So, we have
\begin{equation}
    \tau=0, \qquad G=0, \qquad
    3\tau^*-\tau=0 \quad
    (\text{hence }\quad
    \tau^*=0).
\end{equation}
Therefore, from \eqref{eq:R_B(R)=0}, $M$ is locally
flat.
\hfill$\square$\medskip


\section{ Compact Tricerri-Vanhecke Bochner flat K\"ahler
Surfaces}

 \quad\;
{We give another proof of the result by Kamishima
       \cite {Ka} for the real four dimensional case.}
 Let $M=(M,J,g)$ be a compact Tricerri-Vanhecke Bochner flat
 K\"ahler surface. First, we recall so-called Miyaoka-Yau's
 inequality  \cite {MiYa} :
\begin{equation}
    c_1(M)^2
   \leqq Max\left\{2 \chi(M),3 \chi(M)\right\}.
\end{equation}

\noindent Since $M$ is K\"ahler, the integral formulas \eqref
{eq:p_1(M)} and \eqref{eq:chi(M)} imply
\begin{equation}\label{eq:chi(M)-K}
\begin{aligned}
    \chi(M)
    =&\frac{1}{32\pi^2}\int_M\left\{\|R\|^2-4\|\rho\|^2+\tau^2\right\}dv \\
    =&\frac{1}{32\pi^2}\int_M\left\{\frac{(3\tau^*-\tau)^2}{24}
        -2\Big\|\rho-\frac{\tau}{4}g\Big\|^2
        +\frac{\tau^2}{6}
        +\frac12G\right\}dv\\
    =&\frac{1}{32\pi^2}\int_M\left\{\frac{\tau^2}{3}-2\Big\|\rho-\frac{\tau}{4} g\Big\|^2
    \right\}dv,
\end{aligned}
\end{equation}

\begin{equation}\label{eq:c_1(M)-K}
    c_1(M)^2
   =\frac{1}{32\pi^2}\int_M\left\{\tau^2-4\Big\|\rho-\frac{\tau}{4} g\Big\|^2
   \right\}dv,
\end{equation}
 \noindent respectively. We now assume that $\chi(M) \geqq 0$. Then,
Miyaoka-Yau's inequality implies
\begin{equation}\label{eq:c_1(M)<3-K}
    c_1(M)^2 \leqq 3 \chi(M).
\end{equation}
Then, by \eqref{eq:chi(M)-K}, \eqref{eq:c_1(M)-K}, and
\eqref{eq:c_1(M)<3-K}
 we have
\begin{equation}
   \int_M {2\Big\|\rho-\frac{\tau}{4} g\Big\|^2 }dv\leqq 0
\end{equation}
and hence, $M$ is Einstein (and therefore, in particular, the scalar
curvature of $M$ is constant). Thus, by Theorem \ref{th: K or -K}
and \eqref{eq:R_B(R)=0} we see that $M$ is locally a complex space
form. { Next,} we assume that $\chi(M)<0$. Then, Miyaoka-Yau's
inequality implies
\begin{equation}\label{eq:c_1(M)<2-K}
     c_1(M)^2 \leqq 2 \chi(M).
\end{equation}
Thus, in this case, by \eqref{eq:chi(M)-K},\eqref{eq:c_1(M)-K} and
\eqref{eq:c_1(M)<2-K}, we have
\begin{equation}
   \int_M \frac{\tau^2}{3} ~dv\leqq 0
\end{equation}
and hence $\tau \equiv 0$ on $M$. Thus, by Theorem \ref{th: K or
-K}, we see also that $M$ is locally a product of two oriented
surfaces of constant Gaussian curvatures $K$ and $-K$ ($K \ne 0$).
Summing up the above arguments, we have the following { theorem}.

\begin{thm}
    Let $M=(M,J,g)$ be a compact Tricerri-Vanhecke Bochner flat
    K\"ahler surface.
    Then $M$ is locally a complex space form of complex dimension
    $2$,
    or locally a product of two oriented surfaces
    of different constant Gaussian curvatures $K$ and $-K$
    \rm{(}$K\neq0$\rm{)}.
\end{thm}

\noindent{\bf Remark} {{\rm The above Theorem 5.1 is included in
   the result by Y. Kamishima \cite {Ka} and the proof was first given by
   B. Y. Chen \cite {Chen}. We refer to
\cite{i, Br, Ka} for a further discussion of the Bochner-K\"aher
manifold.
   We may note that our proof is different from theirs.
 }}
\section{ Tricerri-Vanhecke Bochner flat alomst K\"ahler \\Einstein
surfaces}
 \quad\; Let $M=(M, J, g)$ be a compact Tricerri-Vanhecke Bochner flat
almost K\"ahler Einstein surface. From \eqref{eq:R_B(R)=0}, we have
\begin{equation}\label{eq:Eins_R}
\begin{gathered}
    R_{1313}=\frac{3\tau^*-5\tau}{48},\qquad
    R_{1324}=-\frac{3\tau^*-\tau}{48},\qquad
    R_{1414}=\frac{3\tau^*-5\tau}{48},\\
    R_{1423}=\frac{3\tau^*-\tau}{48},\qquad
    R_{1314}=0,\qquad R_{1323}=0.
\end{gathered}
\end{equation}
From \eqref{eq:Eins_R}, we thus have
\begin{equation*}
\begin{gathered}
    u=-R_{1313}+R_{1324}= - \frac{\tau^* - \tau}{8},\\
    v=-R_{1414}-R_{1423}= - \frac{\tau^* - \tau}{8},\\
    w=-R_{1314}-R_{1323}=0,
\end{gathered}
\end{equation*}
and
\begin{equation}
h\equiv (u-v)^2-4w^2=0.
\end{equation}
which implies that $M$ is an almost K\"ahler Einstein surface with
Hermitian Weyl tensor \cite{OS}. Therefore, by virtue of
\cite{ApAm}, { we see immediately that $M$ is a K\"ahler surface.
Therefore, taking account of Theorem \ref{th: K or -K}, we have the
following theorem concerning the Goldberg conjecture \cite{Go,
Se1}}.

\begin{thm}
    Let $M=(M,J,g)$ be a compact Tricerri-Vanhecke Bochner flat
    almost K\"ahler Einstein surface.
    Then $M$ is { locally a complex space form of complex dimension 2}.
\end{thm}

\section{Remarks}
  \quad\;T. Koda \cite{K} has proved that a self-dual
almost Hermitian Einstein surface is a space of pointwise constant
holomorphic sectional curvature. Further, T. Koda and fourth author
of the present paper have proved that a compact self-dual Hermitian
Einstein surface is a complex space form of complex dimension 2
\cite{KoSe}. Therefore, we see that a compact Tricerri-Vanhecke
Bochner flat Hermitian Einstein surface is a complex space form of
complex dimension 2. We herewith introduce an example of a
non-compact Tricerri-Vanhecke Bochner flat Hermitian surface of
pointwise constant holomorphic sectional curvature which is weakly
$*$-Einstein but not Einstein.

\begin{ex}
{\rm    Let $\mathbb{C}$ be the set of complex numbers and $f$ be a
non-constant holomorphic function on
    { $\mathbb{C}^2 = \Bbb C \times \mathbb{C}$}.
    We set $M=\{z=(z_1, z_2) \in \Bbb C^2 \ | \ {\text Re} f(z)>-1 \}$  and assume
    that $M$ is nonempy. Further, we set $u(z)=\text{Re} f(z)$ and  $ \sigma (z) \equiv
    -\text{log}(1+u(z))$.
    Then $\sigma$ is regarded as a smooth function on $M$.
    Let $g$ be the canonical Euclidean metric on ${\mathbb C}^2$ and $J$ be
    the complex structure on $M$ induced by the canonical complex
    structure on ${\mathbb C}^2$. Let $\bar g$ be the Riemannian metric on $M$ defined by
    \begin{equation}
    \bar{g} = e^{2\sigma} g = \frac{1}{(1+u(z))^2}g.
    \end{equation}\\
    Since $M=(M, J, g)$ is a locally flat Hermitian surface (and
    hence, $M$ is a Tricerri-Vanhecke Bochner flat
    Hermitian  surface).
    Since the Tricerri-Vanhecke Bochner curvature tensor $B(R)$ is
    conformally invariant, $(M, J,\bar {g})$ is also Tricerri-Vanheche
    Bochner flat. Further, Tricerri and Vanhecke proved that $(M, J, \bar{g})$ is
    a space of pointwise constant holomorphic sectional curvature
    $c=-e^{2\sigma} \| \text{grad} \sigma \|^2_g$ and $\tau^* = 4c$, where $\| \cdot \|^2_g$
    denotes the square norm with respect to the flat metric $g$ on
    $M$ \cite{GrVa}. We may also check that
Example 4 is a weakly $*$-Einstein manifold.}
\end{ex}

It is known that a Tricerri-Vanhecke Bochner flat almost Hermitian
manifold $M=(M, J, g)$ is a general complex space form if and only
if $M$ is Einstein and weakly $*$-Einstein, and further that a
general complex space form of dimension $2n(\geqq6)$ is locally a
complex space form \cite{TrVa81}. Concerning this result, Lemence
proved that a compact generalized complex space form of dimension
four is locally a complex space form of complex dimension 2
\cite{Le}.


\end{document}